\documentclass[a4paper]{article}
\usepackage[T1]{fontenc}
\usepackage[utf8]{inputenc}
\usepackage{amsmath}
\usepackage{amsthm}
\usepackage{amsfonts}
\usepackage{amssymb}
\usepackage{graphicx}
\usepackage{tikz}
\usepackage{mathtools}
\usepackage[top=3cm, bottom=3cm, left=3cm, right=3cm]{geometry}
\usepackage{calrsfs}
\usepackage{bbm}
\usepackage{algorithm}
\usepackage{cases}
\usepackage[noend]{algpseudocode}
\usepackage{float}
\usepackage{enumerate}
\usepackage{stmaryrd}
\DeclareMathAlphabet{\pazocal}{OMS}{zplm}{m}{n}
\usepackage{lipsum}
\usepackage{enumerate}

\usepackage{mathtools}
\newcounter{mysubequations}

\newcommand{\integerinftyn}{\ \rrbracket \kern-0.20em -\kern-0.20em \infty,n\rrbracket}
\newcommand{\integerninfty}{\llbracket n,+\infty\llbracket}

\newcommand{\integerinftyzero}{\ \rrbracket \kern-0.20em -\kern-0.20em \infty,0\rrbracket}
\newcommand{\integerzeroinfty}{\llbracket 0,+\infty\llbracket}

\newcommand{\integerinftyminusone}{\ \rrbracket \kern-0.20em -\kern-0.20em \infty,-1\rrbracket}

\makeatletter
\newcommand*\bigcdot{\mathpalette\bigcdot@{.5}}
\newcommand*\bigcdot@[2]{\mathbin{\vcenter{\hbox{\scalebox{#2}{$\m@th#1\bullet$}}}}}
\makeatother

\newcommand{\mydef}{\vcentcolon=}

\newtheorem{theorem}{\indent Theorem}[section]
\newtheorem{proposition}[theorem]{\indent Proposition}

\newtheorem{remark}[theorem]{\indent Remark}
\newtheorem{lemma}[theorem]{\indent Lemma}

\newtheorem{definition-theorem}[theorem]{\indent Definition-Theorem}

\makeatletter
\renewenvironment{proof}[1][\proofname]{%
	\par\pushQED{\qed}\normalfont%
	\topsep6\p@\@plus6\p@\relax
	\trivlist\item[\hskip\labelsep\bfseries#1\@addpunct{.}]%
	\ignorespaces
}{%
	\popQED\endtrivlist\@endpefalse
}
\makeatother

\def \N{\mathbb{N}}
\def \R{\mathbb{R}}
\def \Z{\mathbb{Z}}
\def \V{\mathbb{V}}

\def \G{\mathcal{G}}

\def \L{\mathcal{L}}

\def \P{\mathbb{P}}
\newcommand{\E}{\mathbb{E}}
\newcommand{\Var}{\mathbb{V}\mathrm{ar}}
\newcommand{\Cov}{\mathbb{C}\mathrm{ov}}

\title{\bf Convergence of the Temporal Averages of a Metastable System of Spiking Neurons}

\author{Morgan André \\ \\ \textit{IMECC-UNICAMP}}

\begin{document}

\maketitle

\begin{abstract}
We consider a stochastic system of spiking neurons which was previously proven to present a metastable behavior for a suitable choice of the parameter, in the sense that the time of extinction is asymptotically memory-less when the number of components in the system goes to $\infty$. In the present article we complete this work by showing that, previous to extinction, the system tends to stabilize in the sense that temporal means taken on an appropriate time scale converge in probability to some fixed value. This property is sometime called \textit{thermalization}.
\end{abstract}

\vspace{0.4 cm}



\vspace{1 cm}

\section{Introduction}

Informally the model is as follows. Each element in the system (a neuron) is associated to a random variable called its membrane potential. Each neuron is also associated to a point process which intensity varies across time depending on the current value of the membrane potentials, representing the spiking times, and to another point process (which is Poisson of some fixed rate $\gamma \geq 0$) representing the leakage phenomenon, that is, the drift of the membrane potential toward its resting value caused by the natural diffusion of ions through de membrane when some equilibrium has not been reached. Thus, unlike in the original model the present one was inspired from, these leaks occurs at discrete times at which the membrane potential is reset to its resting value (conventionally set to $0$). Moreover, like in the original model, the membrane potential of any given neuron is also reset to the resting value whenever this neuron spikes. Furthermore all the neurons are excitatory, with the same synaptic weight (conventionally set to $1$). We give a formal definition of the system in the next section.

This model was proposed in \cite{ferrari18}, as one of the continuous time versions of the model introduced in \cite{galves2013}, sometimes called the Galves-Löcherbach model, or simply GL model. It was studied for the infinite one-dimensional lattice instantiation in which the neurons are indexed by $\Z$ and each neuron is connected to its immediate neighbors, on the right and on the left, and was proven to present a phase transition with respect to the parameter $\gamma$: there exists a critical threshold for this parameter such that if $\gamma$ is taken above it then each neuron stop emitting spikes in finite time almost surely, while if $\gamma$ is taken below the threshold then each neuron has a positive probability of emitting spikes forever.

Then the study of the metastable properties of this model was initiated in \cite{andre19}, in an attempt to fill the gap for a mathematical treatment of the subject of metastability in neuroscience, which is central but rarely treated from the rigorous and microscopical perspective of statistical mechanics (on the role of the concept of metastability in neuroscience, which is out of the scope of this article, see for example \cite{werner2007}). It was proven there that in at least a sub-region of the sub-critical region, if we consider the finite version of the model described above with $n$ neurons (i.e. the neurons are indexed on a finite windows such as $\Z \cap [ -n , n]$) then the instant of the last spike of the system (which is almost surely finite) converges to an exponential random variable as $n$ diverges. It is the first of the two characteristic properties of metastable dynamics as formalized in the seminal paper \cite{cassandro1984}. Then this type of convergence was also proved to hold for other versions of the GL model  (see \cite{planche} and \cite{locherbach}). In the present article the result of \cite{andre19} is slightly improved in Section \ref{sec:mainthm}, as a side consequence of the analysis that will be conducted up to that point is that the convergence toward an exponential actually holds in the whole sub-critical region.

Our main result (Theorem \ref{thm:main}) complete the study of the metastable properties of this system by showing that, before extinction, the spiking activity stabilizes in a dynamic which resemble stationarity. This is the second of the two properties characterizing metastability, sometime called \textit{thermalization} (see \cite{schonmann85}). More precisely we show that in the finite version mentioned above the temporal averages of spikes along time, computed on a suitable time scale, converge in probability as $n$ diverges to the asymptotic spatial average of active neurons in the infinite system. In that sense it can also be seen as a result of \textit{ergodicity}. One of the important ideas of the proof of the main result is to consider an auxiliary process, namely the spiking rate process, which is an interacting particle system taking value in $\{0,1\}^\Z$, as well as its \textit{dual} (in the sense introduced by T. Harris in \cite{harris76}). Then exploiting interesting properties of this auxiliary process allows us to derive the proof of our central result.

The paper is organized as follows. In Section \ref{sec:themodel} we introduce formally our model and state the main result. In Section \ref{sec:auxdual} we introduce the auxiliary process as well as its dual; we also list previously obtained results concerning these processes which will be important in the sequel. In Section \ref{sec:positivedrift} we obtain a result about the drift of the right-most component of the dual process. This result is then used in Section \ref{sec:subexpest} to prove that the auxiliary process has exponentially decaying time correlations, which is a crucial ingredient of the proof of the main result. Finally this main result is proven in Section \ref{sec:mainthm}.

\section{Definition of the model and main result}
\label{sec:themodel}

The general model considered in this work is as follows. Let $I$ be a finite or countable set representing the neurons, and to each $i \in I$ associate a set $\V_i \subset S$ of \textbf{presynaptic neurons}. Each neuron $i \in I$ has a \textbf{membrane potential} evolving over time, represented by a stochastic process which takes its values in the set $\mathbb{N}$ of non-negative integers and which is denoted $(X_i(t))_{t \geq 0}$. The evolution of $(X_i(t))_{t \geq 0}$ depends on two type of point processes, denoted $(N^{\dagger}_i(t))_{t \geq 0}$ and $(N_i(t))_{t \geq 0}$ for each neuron $i \in I$. $(N^{\dagger}_i(t))_{t \geq 0}$ is a Poisson process of some parameter $\gamma$, representing the \textbf{leak times}. At any of these leak times the membrane potential of the neuron concerned is reset to $0$. Finally $(N_i(t))_{t \geq 0}$ represents the \textbf{spiking times}, and its infinitesimal rate at time $t$ is given by $\phi(X_i(t))$, where $\phi$ is some rate function. When a neuron spikes its membrane potential is reset to $0$ and the membrane potentials of all of its \textbf{postsynaptic neurons} (that is the neurons of the set $\{ j : i \in  \mathbb{V}_j\}$) are increased by one. All the point processes involved are assumed to be mutually independent.

\vspace{0.4 cm}

Mathematically, beside asking that $(N^{\dagger}_i(t))_{t \geq 0}$ be a Poisson process of some parameter $\gamma$, this is the same as saying that $(N_i(t))_{t \geq 0}$ is the point process characterized by the two following equations $$ \label{spikephi} \E (N_i (t)- N_i (s)|\mathcal{F}_s) = \int_s^t \E (\phi (X_i(u))|\mathcal{F}_s)du $$ where $$ X_i(t) = \sum_{j \in \mathbb{V}_i}\int_{]L_i(t),t[}dN_j(s),$$ $L_i(t)$ being the time of the last event affecting neuron $i$ before time $t$, that is, $$ L_i(t) = \sup \Big\{s \leq t : N_i(\{s\}) = 1 \text{ or } N^\dagger_i(\{s\}) = 1\Big\}.$$ $(\mathcal{F}_t)_{t\geq 0}$ is the standard filtration generated by the family $\{ N_i(s), N_i^\dagger(s), s \leq t, i \in S\}$.

\vspace{0.4 cm}

We continue the study initiated in \cite{ferrari18} and \cite{andre19}, and we study the specific case in which the activation function is simply a hard threshold of the form $\phi(x) = \mathbbm{1}_{x>0}$ and the spatial structure of the network is given by a nearest-neighbor interaction on the one-dimensional lattice, that is we set $I=\Z$ and $\mathbb{V}_i = \{i-1,i+1\}$. As stated in the introduction, it was proven in \cite{ferrari18} that there exists some critical value $0 < \gamma_c < \infty$ such that, assuming that the system start in a state in which every neuron has a positive membrane potential, then for any $i \in I$

$$\P \Big( N_i([0,\infty[) \text{ } < \infty \Big) = 1 \text{ if } \gamma > \gamma_c,$$

while

$$ \P \Big( N_i([0,\infty[) \text{ } = \infty \Big) > 0 \text{ if } \gamma < \gamma_c.$$

\vspace{0.4 cm}

This infinite system is interesting in itself but our main concern in this article (as in \cite{andre19}) is its finite counterpart, that is the system in which the neurons are indexed on $I_n = \llbracket -n , n \rrbracket$ for some $n \in \N$ (Here and in the sequel $\llbracket -n,n \rrbracket$ is a short-hand for $\Z \cap [-n,n]$) and in which the sets of presynaptic neurons for each neuron $i \in I_n$ is given by $\V_{n,i} = \{-i,i\} \cap \llbracket -n , n \rrbracket$. The reason for introducing the infinite system and actually studying some of its properties (Sections \ref{sec:positivedrift} and \ref{sec:subexpest}) is that these are extremely useful in order to obtain the finite-case result we are interested in. What we prove is that in this finite system, in the sub-critical regime, and if $n$ is big enough, then counting the number of spikes occurring in a given time interval before extinction, for a given subset of neurons of interest in the system, shall give a number which with high probability is close to some fixed value, which depends only on the parameter $\gamma$ (and of course of the number of neurons in the subset considered). This property captures the pseudo-stationarity which is characteristic of metastable systems. 

\vspace{0.4 cm}

More precisely, let $F \subset \Z$ with $|F| < \infty$. Then we define, for any $t, R \in \R_+$ and for any $n \geq 0$, the following quantity $$ \widehat{N}^n_R \left(t,F\right) = \frac{1}{R} \sum_{i \in F \cap I_n} N_i \left( [t, t+R] \right),$$ where the superscript $n$ here indicates that we are considering the finite system defined with respect to the $I_n$ and $\V_{n,i}$ defined above. $\widehat{N}^n_R \left(t,F\right)$ is the average number of spikes emitted by the neurons in $F$ on a time interval of length $R$, starting the enumeration a time $t$. Our main result is the following theorem.

\vspace{0.4 cm}

\begin{theorem} \label{thm:main}
	Suppose $0 < \gamma < \gamma_c$ and let $(R_n)_{n \geq 0}$ be an increasing sequence of positive real numbers satisfying $$ R_n \underset{n \rightarrow \infty}{\longrightarrow} + \infty \ \ \text{ and } \ \ \frac{R_n}{\E \left( \tau_n \right)} \underset{n \rightarrow \infty}{\longrightarrow} 0.$$
	
	There exists some $0 < \rho < 1$ (which depends only on $\gamma$) such that for any $t \geq 0$ $$ \widehat{N}^n_{R_n} \left(t,F\right) \overset{\P}{\underset{n \rightarrow \infty}{\longrightarrow}} |F| \cdot \rho.$$
\end{theorem}

\vspace{0.4 cm}

\section{The auxiliary process and its dual}

\label{sec:auxdual}

Consider the infinite process introduced above, with $I = \Z$ and $\V_i = \{i-1,i+1\}$. For any $i \in \Z$ and $t \geq 0$ we write $\xi_i(t) = \mathbbm{1}_{X_i(t) > 0}$, and $\xi(t) = (\xi_j(t))_{j \in \Z}$. The resulting \textbf{auxiliary process} $(\xi(t))_{t \geq 0}$ gives the state of each neuron at any time among the two possibilities: \textbf{active} or susceptible to spike if $\xi_i(t) = 1$, and \textbf{quiescent} or not susceptible to spike if $\xi_i(t) = 0$. This stochastic process is an \textit{interacting particle systems} (see \cite{liggett85}) and it's dynamic can be described as follows: any active neuron spikes at rate $1$, and when is does so it "activates" its two neighbors (in case they weren't already) while it immediately turns itself quiescent, moreover it also spontaneously becomes quiescent at rate $\gamma$ because of the leakage. In a more formal way $(\xi(t))_{t \geq 0}$ is the Markovian process on $\{0,1\}^\Z$ which infinitesimal generator is given by

\begin{equation} \label{generator} \L f (\eta) = \gamma \sum_{i \in \Z} \Big(f(\pi^{\dagger}_i(\eta)) - f(\eta)\Big) + \sum_{i \in \Z} \eta_i\Big(f(\pi_i(\eta)) - f(\eta)\Big),
\end{equation} where $f : \{0,1\}^{\Z} \rightarrow \R$ is a cylinder function\footnote{Here and in the rest of the paper we call \textit{cylinder function} any function $f: \{0,1\}^{\Z} \rightarrow \R$ which only depends on a finite number of sites. The set $S \subset \Z$ of the sites on which $f$ depends is called the \textit{support} of $f$.} and the $\pi^{\dagger}_i$'s and $\pi_i$'s are maps from $\{0,1\}^{\Z}$ to $\{0,1\}^{\Z}$ defined for any $i \in \Z$ as follows:

\begin{equation} \label{leaksmap}
    {\Big(\pi^{\dagger}_i(\eta)}\Big)_j= 
\begin{cases}
    0& \text{if } j = i,\\
    \eta_j & \text{otherwise},
\end{cases}
\end{equation} and

\begin{equation} \label{spikesmap}
    {\Big(\pi_i(\eta)}\Big)_j= 
\begin{cases}
    0 & \text{if } j = i,\\
    \max (\eta_i, \eta_j)  & \text{if } j \in \{i-1, i+1\},\\
    \eta_j & \text{otherwise}.
\end{cases}
\end{equation}

\vspace{0.4 cm}

\subsection{Graphical construction}

\label{subsec:graphconsaux}

It is possible to use a graphical construction of the type of the construction introduced by Harris in \cite{harris78} to propose an alternative definition of the auxiliary process. The construction is as follows. For any neuron $i \in \Z$, we let $(N^*_i(t))_{t \geq 0}$ and $(N^\dagger_i(t))_{t \geq 0}$ be two independent homogeneous Poisson processes with respective intensity $1$ and $\gamma$. Let $(T^*_{i,n})_{n \geq 0}$ and $(T^\dagger_{i,n})_{n \geq 0}$ be their respective jump times. All the Poisson processes are assumed to be mutually independent and we let $(\Omega, \mathcal{F}, \P)$ be the probability space on which these Poisson processes are defined.

\vspace{0.4 cm}

Now we adjoin the following structure to the time-space diagram $\Z \times \R_+$:

\begin{itemize}
\item for all $i \in \Z$ and $n \in \N$ put a `$\delta$' mark at the point $(i,T^{\dagger}_{i,n})$,
\item for all $i \in \Z$ and $n \in \N$ put an arrow pointing from $(i,T^*_{i,n})$ to $(i+1,T^*_{i,n})$ and another pointing from $(i,T^*_{i,n})$ to $(i-1,T^*_{i,n})$.
\end{itemize}

\vspace{0.4 cm}
We obtain a random structure that we denote $\G$, which consists of the time-space diagram $\Z \times \R$ augmented by the set of `$\delta$' marks and horizontal arrows. Then for any $i,j \in \Z$ and $t < s$ we call a \textbf{path} from $(i,t)$ to $(j,s)$ on $\G$ any alternated sequence of contiguous closed time segments and arrows, starting and ending by a time segment (possibly reduced to a single point), such that $(i,t)$ is the bottom endpoint of the first segment, and $(j,s)$ the top endpoint of the last segment. Moreover we say that a path is \textbf{valid} if 

\begin{enumerate}[(i)]
\item none of the time segments contains a `$\delta$' mark,
\item none of the time segments contains a $(i,T^*_{i,n})$ point in its interior or at its bottom endpoint.
\end{enumerate}

In other words a path is deemed valid if, as you goes upward along this path, you never encounter a `$\delta$' mark neither you cross the rear side of an arrow. For any two points $i$ and $j$ in $\Z$ we write $(i,t) \longrightarrow (j,s)$ if there is a valid path from $(i,t)$ to $(j,s)$. For any two sets $A,B \subset \Z$ we also write $A \times t \longrightarrow  (j,s)$ (resp. $(i,t) \longrightarrow  B \times s$) if there exists some $i$ in $A$ (resp. some $j$ in $B$) such that $(i,t) \longrightarrow  (j,s)$, and we write $A \times t \longrightarrow  B \times s$ when $(i,t) \longrightarrow  (j,s)$ for $i \in A$ and $j \in B$.

\vspace{0.4 cm}

With this construction we can easily give the following characterization of our stochastic process\footnote{Notice that with this new definition the process is defined on $\mathcal{P}(\Z)$ instead of $\{0,1\}^{\Z}$. It is of course only a matter of notation, as any element $\eta$ of $\{0,1\}^{\Z}$ can be bijectively mapped to an element $A$ of $\mathcal{P}(\Z)$ via the relation $A = \{i \in \Z \text{ such that } \eta_i = 1\}$. In practice we will indifferently use both ways.}. For any $A \in \mathcal{P}(\Z)$, and for any $t \geq 0$ :

$$\xi^A(t) = \{j \in \Z :  A \times 0 \longrightarrow (j,t)\}.$$

\vspace{0.4 cm}

Then $(\xi^A(t))_{t \geq 0}$ is the process with generator  (\ref{generator}) and initial state $\xi^A(0) = A$.

\vspace{0.4 cm}

Notice that once all the neurons are quiescent the system will remains quiescent for eternity. In other words $\emptyset$ is an absorbing state. Therefore we can define the \textbf{extinction time} of the system. For any $A \in \mathcal{P}(\Z)$ the time of extinction, denoted $\tau^A$, is defined as follows $$\tau^A = \inf \{t \geq 0: \xi^A(t) = \emptyset\}.$$

\vspace{0.4 cm}

\subsection{Dual process} 
\label{sec:dualprocess}
It is possible to define a \textbf{dual process} for $(\xi(t))_{t \geq 0}$, which is particularly useful for the study of the original process. Again, for any $i \in \Z$, let's consider two independent homogeneous Poisson processes $(\tilde{N}^*_i(t))_{t \geq 0}$ and $(\tilde{N}^\dagger_i(t))_{t \geq 0}$ with intensity $1$ and $\gamma$ respectively, and let $(\tilde{T}^*_{i,n})_{n \geq 0}$ and $(\tilde{T}^\dagger_{i,n})_{n \geq 0}$ be their respective jump times. As previously all the Poisson processes are assumed to be mutually independent.

\vspace{0.4 cm}

The time-space diagram $\Z \times \R_+$ is then augmented in order to obtain the dual structure $\tilde{\G}$ as follows:

\begin{itemize}
\item for all $i \in \Z$ and $n \in \N$ put a `$\delta$' mark at the point $(i,\tilde{T}^{\dagger}_{i,n})$,
\item for all $i \in \Z$ and $n \in \N$ put an arrow pointing from $(i+1,\tilde{T}^*_{i,n})$ to $(i,\tilde{T}^*_{i,n})$ and another pointing from $(i-1,\tilde{T}^*_{i,n})$ to $(i,\tilde{T}^*_{i,n})$.
\end{itemize}

Now we say that a path in $\tilde{\G}$ is a \textbf{dual-valid path} if it satisfies the following constraints:

\begin{enumerate}[(i)]
	\item none of the time segments contains a `$\delta$' mark,
	\item none of the time segments contains a $(i,\tilde{T}^*_{i,n})$ point in its interior or at its top endpoint.
\end{enumerate}

In other words a path is dual-valid if, as you goes upward along this path, you never encounter a `$\delta$' mark neither you cross the tip of an arrow. For any two points $i$ and $j$ in $\Z$ we write $(i,t) \overset{\text{dual}}{\longrightarrow} (j,s)$ if there is a dual-valid path from $(i,t)$ to $(j,s)$, and we adopt similar notations when sets are considered instead of points.

\vspace{0.4 cm}

Then for any $A \in \mathcal{P}(\Z)$ and for any $t \geq 0$ we write $$\eta^A(t) = \{j \in \Z :  A \times 0 \overset{\text{dual}}{\longrightarrow} (j,t)\}.$$ 

\vspace{0.4 cm}

The process $(\eta(t))_{t \geq 0}$ thus defined is the dual of $(\xi(t))_{t \geq 0}$. Its dynamic can be briefly described as follows: an active neuron "activates" its neighbors at rate $1$ (on the left and on the right independently), while it turns itself quiescent at rate $\gamma$ or $1+\gamma$ depending on whether it has active neighbors or not. 

\vspace{0.4 cm}

We briefly explain the point of introducing a dual process, and we refer to \cite{harris76} and \cite{bertein} for more details about duality. The crucial point is that we can relate the initial process $(\xi(t))_{t \geq 0}$ and its dual in the following way. We fix some $s \in \R^+$, and for any $0 \leq t \leq s$ and $A \in \mathcal{P}(\Z)$ we define the following random variable on $\mathcal{P}(\Z)$ via the random graph $\mathcal{G}$ of the previous section 

\begin{equation} \label{eq:defdual}
\zeta^A(t) = \{ i \in \Z: (i,s-t) \longrightarrow A \times s\}.
\end{equation}

That way it is easy to see that $(\zeta(t))_{t \in [0,s]}$ is the dual process $(\eta(t))_{t \geq 0}$ restricted to the time interval $[0,s]$, but built on the probability space $(\Omega, \mathcal{F}, \P)$ of the orginal process $(\xi(t))_{t \geq 0}$. We call this way of defining the dual process by a coupling with the initial process---with the time reversed---the \textbf{backward version} of the dual process. Moreover it is straightforward to check that for any $A,B \in \mathcal{P}(\Z)$ and any $t \geq 0$ the following holds
\begin{equation} \label{eq:dualityrelation}
\{ \xi^A(s) \cap B \neq \emptyset\} = \{ \zeta^B(s) \cap A \neq \emptyset \}.
\end{equation}

Thus, the following proposition holds.

\begin{proposition} \label{prop:duality}
For any $A, B \in \mathcal{P}(\Z)$, and $t \geq 0$ we have

$$ \P \Big( \xi^A(t) \cap B \neq \emptyset \Big) = \P \Big( \eta^B(t) \cap A \neq \emptyset \Big).$$
\end{proposition}

\vspace{0.4 cm}

For any $A \in \mathcal{P}(\Z)$ we denote by $\sigma^A$ the time of extinction of the dual process $$\sigma^A = \inf \{t \geq 0: \eta^A(t) = \emptyset\}.$$ This time of extinction will be of great importance, as an important ingredient of the proof of the main theorem consists in using the exponential decay of the time correlations of the original interacting particle system (Theorem \ref{th:timecorrelation}), which is proven using the duality property and an exponential bound on this extinction time (Theorem \ref{th:exptinf}).  

\subsection{Important properties of the auxiliary process and its dual}
\label{subsec:importantprop}

In this section we summarize some of the properties that have already been proven for the interacting particle system and its dual. In order to avoid redundancy most of the properties are stated only for the process $(\xi(t))_{t \geq 0}$ but they hold for the dual process $(\eta(t))_{t \geq 0}$ as well. Moreover the phase transition property is stated in terms of the extinction time of $(\eta(t))_{t \geq 0}$. The proofs can be found in Section 4 of \cite{andre19}. Below, and from now on, we will use some new notation. The time will sometimes be written as a subscript if it is more suitable, writing $(\xi_t)_{t \geq 0}$ instead of $(\xi(t))_{t \geq 0}$. We sometimes write $\xi \equiv 0$ (or $\eta \equiv 0$) for the state in which all neurons are quiescent, and $\xi \equiv 1$ (or $\eta \equiv 1$) for the state in which all neurons are active. When the initial state is a singleton we drop the curly bracket, writing for example $\sigma^0$ instead of $\sigma^{\{0\}}$. When the initial state is the whole lattice $\Z$ we will omit the superscript, writing simply $\xi(t)$ for $\xi^\Z(t)$. Moreover the state space $\{0,1\}^\Z$ is associated with the partial order relation defined for any $\xi, \eta \in \{0,1\}^\Z$ by: $\xi \leq \eta$ if and only if $\xi_i \leq \eta_i$ for all $i \in \Z$. Whenever we say that a function on $\{0,1\}^\Z$ is monotonous, it is to be understood with respect to this partial order. Finally, for any probability measure $\nu$ on $\{0,1\}^\Z$ (associated with its standard Borel $\sigma$-algebra) and for any measurable function $f$ on $\{0,1\}^\Z$ we write $\nu(f) = \int f d\nu$.

\begin{enumerate}[(i)]
    \item \textbf{Additivity}: From the graphical constructions we immediately obtain that for any $A,B \in \mathcal{P}(\Z)$ and for any $t \geq 0$ the following holds
    \begin{equation} \label{additivity}
    \xi^{A \cup B}(t) = \xi^A(t) \cup \xi^B(t).
    \end{equation}
    
    \item \textbf{Monotonicity}: The previous property implies that for any $A,B \in \mathcal{P}(\Z)$ such that $A \subset B$ and for any $t \geq 0$
    \begin{equation} \label{monotonicity}
    \xi^{A}(t) \subset \xi^B(t).
    \end{equation}
    
    \item \textbf{Attractiveness}: By definition an interacting particle system on $\{0,1\}^\Z$ with semi-group $(S(t))_{t \geq 0}$ is attractive if for any increasing function on $\{0,1\}^\Z$ the function $S(t)f$ is increasing for any $t \geq 0$. For any  $\xi, \eta \in \{0,1\}^\Z$ satisfying $\xi \leq \eta$ it is immediate using monotonicity that for any increasing $f$ and for any $t \geq 0$ we have $\E^\xi\left(f(\xi_t)\right) \leq \E^\eta\left(f(\xi_t)\right)$, so that our system is indeed attractive.
    
    \item \textbf{Translation invariance}: It is clear from the graphical construction that the law of the process does not change if the time-space diagram is translated to the right or to the left.
    
    \item \textbf{Phase transition}: The phase transition stated in Section \ref{sec:themodel} can be restated in term of the dual process:
    
    $$\P \left(\sigma^0 = +\infty \right) > 0 \text{ if } \gamma < \gamma_c$$
    and 
    $$\P \left(\sigma^0 = +\infty \right) = 0 \text{ if } \gamma > \gamma_c.$$
    
    \item \textbf{Invariant measures}: If $\gamma < \gamma_c$, then there exists a non-trivial invariant measure (in the sense that it doesn't give mass $1$ to $\xi \equiv 0$) for $(\xi_t)_{t \geq 0}$, which corresponds to the weak limit of $\xi_t$ when $t$ diverges, and which we denote $\mu$. There is an analogous invariant measure for the dual process $(\eta_t)_{t \geq 0}$ which we denote $\tilde{\mu}$.
    
    \item \textbf{Stochastic monotonicity}: The convergence toward $\mu$ is monotonous in the sense that for any continuous and increasing function $f:\{0,1\}^\Z \rightarrow \R$ and for $0 \leq s < t$ the following holds $$ \E\left( f(\xi_s)\right) \geq \E\left( f(\xi_t)\right) \geq \mu(f).$$ 
    
    \item \textbf{Positive density}: Define the density of the system\footnote{Notice that, by translation invariance, the fact that we define $\rho$ with respect to neuron $0$ is purely conventional.} $\rho = \mu \left( \{\eta : \eta_0 = 1 \} \right)$. By phase transition and duality (Proposition \ref{prop:duality}), if $\gamma < \gamma_c$ then $\rho > 0$. The same holds for the density of the dual process $\tilde{\rho} = \tilde{\mu} \left( \{\eta : \eta_0 = 1 \} \right)$. While this result isn't proven for the dual in \cite{andre19}, as it wasn't explicitly  needed, it is very easy to prove it using other results proven there, as shown by the following computation. Let $H \subset \{0,1\}^\Z$ and suppose $\gamma < \gamma_c$, then by Proposition 4.8 from \cite{andre19} $$\tilde{\mu} (H) \leq \tilde{\mu} (\eta \equiv 0) + \sum_{i \in \Z} \tilde{\mu} \left( \{\eta \in H: \eta_i = 1\}\right) \leq \sum_{i \in \Z} \tilde{\rho},$$ so that if $\tilde{\rho} = 0$ then $\tilde{\mu}$ shall be identically equal to $0$, which is obviously a contradiction (to the fact that it shall be a probability measure for example).
    
    \item \textbf{Spatial ergodicity}: The measure $\mu$ is spatially ergodic. See Theorem 4.9 in \cite{andre19} for more details. While this Theorem is proven only for $\mu$ there, it is easy to check that all the arguments hold for $\tilde{\mu}$ as well.
 
\end{enumerate}

A fact that will be important in the two next sections is that $(\xi(t))_{t \geq 0}$ and $(\eta(t))_{t \geq 0}$ fall into the category of what is called \textbf{growth models} in \cite{durrett80}, that is, attractive and translation invariant systems with $\emptyset$ as an absorbing state and finite range interaction.

\section{Preliminary results}
\label{sec:positivedrift}
In this section we study the drift of the edge of the dual process $(\eta(t))_{t \geq 0}$. The main result is that in the sub-critical regime the drift is linear, with a positive slope (Proposition \ref{prop:convtoalpha} and Proposition \ref{prop:posalpha}). This fact will be of importance in order to prove exponential estimates for the time of extinction and the time correlations.

\vspace{0.4 cm}

For any set $A \in \mathcal{P}(\Z)$ we write: $$r^A_t = \max \left\{i\in\eta^A_t\right\} \ \text{ and } \ l^A_t = \min\left\{i\in\eta^A_t\right\}.$$

Moreover we write $(\eta^-_t)_{t \geq 0}$ and $(\eta^+_t)_{t \geq 0}$ for the dual processes starting from $\eta^-_0 = \integerinftyzero$ and $\eta^+_0 = \integerzeroinfty$ respectively, and for any $t \geq 0$ we denote $r^-_t = \max \left\{i\in\eta^-_t\right\}$ and $l^+_t = \min\left\{i\in\eta^+_t\right\}.$

\vspace{0.4 cm}

Let's start with the following lemma.

\begin{lemma} \label{lemma:threelemmanotdead}
For any $t \geq 0$, if $\sigma^0 > t$ then:

\begin{enumerate}[(i)]
    \item $r^0_t = r^-_t$ and $l^0_t = l^+_t$,
    \item $\eta^0_t = \eta^-_t \cap [ l^0_t, r^0_t ] = \eta^+_t \cap [ l^0_t, r^0_t ]$,
    \item and $l^+_s \leq r^-_s$  for every $s \leq t$.
\end{enumerate}
\end{lemma}

\begin{proof}
The proof of $(i)$ follows easily from the graphical construction and is quite similar to the proofs of the three first lemmas in Section 5 of \cite{andre19}. We give a quick sketch here. By set monotonicity $r^0_t \leq r^-_t$. Moreover, as $r^-_t \in \eta^-_t$, by definition we can find a dual-valid path from some $i \in \integerinftyzero$ to $r^-_t$. If $i=0$ then $r^-_t \in \eta^0_t$ so that $r^0_t \geq r^-_t$ and the proof is over. Now if $i<0$ then the dual valid-path from $i$ to $r^-_t$ has to cross the left border of $(\eta^0_s)_{s \geq 0}$ somewhere before $t$. Now the concatenation of the path following the left border from $0$ to the crossing point and the path going from the crossing point to $r^-_t$ is a dual-valid path from $0$ to $r^-_t$, proving again that $r^0_t \geq r^-_t$, which ends the proof. Of course the fact that $l^0_t = l^+_t$ is then immediate by symmetry.

\vspace{0.4 cm}

For $(ii)$ by symmetry we only need to prove $\eta^0_t = \eta^-_t \cap [ l^0_t, r^0_t ]$, the proof for the other equality being obviously identical. The fact that $\eta^0_t \subset \eta^-_t \cap [ l^0_t, r^0_t ]$ is immediate by monotonicity so that it suffices to show the reverse inclusion. This is easily done by the same kind of argument as in the first item. Let $j \in \eta^-_t \cap [ l^0_t, r^0_t ]$, then there exists $i \in \integerinftyzero$ such that $(i,0) \overset{\text{dual}}{\longrightarrow} (j,t)$. If $i=0$ then there is nothing to prove, while if $i<0$ then the path from $i$ to $j$ has to cross $\{l^0_t: t \geq 0\}$ somewhere implying that $j \in \eta^0_t$.

\vspace{0.4 cm}

For the last item notice that if $\eta^0_t \neq \emptyset$ then $\eta^0_s \neq \emptyset$ for every $s \leq t$ so that by $(i)$, for any fixed $s \leq t$ $$ l^+_s = l^0_s \leq r^0_s = r^-_s.$$  
\end{proof}

\begin{remark}
It is easy to see that in the above lemma one can replace the initial state $\{0\}$ by any set $A \subset \Z^-$ such that $0 \in A$ (resp. $A \subset \Z^+$ such that $0 \in A$) and that item $(i)$ and $(ii)$ still hold.
\end{remark}

\vspace{0.4 cm}

We also have the following useful lemma.

\begin{lemma} \label{lemma:diffsup1}
For any infinite set $A \subset \Z$ having finitely many positive elements and some $i \in \Z$ satisfying $i > \max\{j \in A\}$, the following holds for all $t \geq 0$ $$\E\left( r^{A\cup\{i\}}_t - r^A_t\right) \geq 1.$$
\end{lemma}

\begin{proof}
It is a direct consequence of the additivity of our process as well as its translation invariance. See the proof of Lemma 2.21 on page 282 of \cite{liggett85}.
\end{proof}

\vspace{0.4 cm}

We have the following

\begin{proposition} \label{prop:convtoalpha}
There exists a constant $\alpha(\gamma) \in [-\infty,\infty[$ such that the following holds $$ \frac{r^-_t}{t} \underset{t \rightarrow \infty}{\longrightarrow} \alpha(\gamma) \text{ almost surely.}$$

Moreover, if $\gamma < \gamma_c$ then $\alpha(\gamma) \geq 0$ and the convergence occurs in $L^1$.
\end{proposition}

\begin{proof}
The existence of $\alpha(\gamma)$ and the almost sure convergence follow from Theorem 2.1 in \cite{durrett80}, since our process is a growth model. Furthermore, if $\alpha(\gamma) < 0$ then almost surely $r^-_t \underset{t \rightarrow \infty}{\longrightarrow} -\infty$ and by symmetry $l^+_t \underset{t \rightarrow \infty}{\longrightarrow} +\infty$ so that using item $(iii)$ in Lemma \ref{lemma:threelemmanotdead} we have 

\begin{equation} \label{eq:tozero}
\P\big(\eta^0_t \neq \emptyset \big) \underset{t \rightarrow \infty}{\longrightarrow} 0.
\end{equation}

But by duality we know that 
\begin{equation} \label{eq:torho}
\P\big(\eta^0_t \neq \emptyset \big) = \P\big(\xi_t (0) = 1\big)  \underset{t \rightarrow \infty}{\longrightarrow} \rho_\gamma,
\end{equation}

and as $\rho_\gamma > 0$ in the sub-critical regime (\ref{eq:tozero}) and (\ref{eq:torho}) implies that $\alpha(\gamma) \geq 0$ when $\gamma < \gamma_c$. Then the $L^1$ convergence follows from the second part of Theorem 2.1 in \cite{durrett80}.
\end{proof}

\vspace{0.4 cm}

The following lemma, which uses Lemma \ref{lemma:diffsup1}, states that in the sub-critical regime the drift coefficient decreases with $\gamma$ in a superlinear fashion.

\begin{lemma} \label{lemma:superlin}
For any $\gamma$ and $\lambda$ such that $0 \leq \gamma + \lambda < \gamma_c$ we have $$\alpha(\gamma) - \alpha(\gamma + \lambda) \geq \lambda.$$
\end{lemma}

\begin{proof}
We adapt an argument from \cite{durrett80}. We build different versions of the dual process on the same probability space using the graphical construction. Suppose that for each integer $k \in \Z$ we have a Poisson process of intensity $\gamma$ and a Poisson process of intensity $\lambda$. We use both of these processes to put `$\delta$' marks on the time-space diagram. Let $(\eta_t(\gamma))_{t \geq 0}$ denote the process which use only the marks coming from the $\gamma$ Poisson processes, and let $(\eta_t(\gamma + \lambda))_{t \geq 0}$ be the process which use the marks of both family of Poisson processes. We let $(r^-_t(\gamma))_{t \geq 0}$ and $(r^-_t(\gamma+\lambda))_{t \geq 0}$ denote the corresponding right edge processes. Then we define the following stopping time:

$$\kappa = \inf\{t \geq 0: r^-_t(\gamma+\lambda) < r^-_t(\gamma)\}.$$

We define a third process, denoted $(\widehat{\eta}_t)_{t \geq 0}$, which use only marks coming from the Poisson processes of parameter $\gamma$ up to time $\kappa$, and then the marks coming from both family of Poisson processes. We require that the initial state of this process is $\widehat{\eta}_0 = \integerinftyzero$. Moreover let $(\widehat{r}_t)_{t \geq 0}$ denotes its right edge. For any $t \geq 0$ we have $$ \E \left(r^-_t(\gamma) - r^-_t(\gamma+\lambda)\right) \geq \E \left(\widehat{r}_t - r^-_t(\gamma+\lambda)\right) \geq \E \left((\widehat{r}_t - r^-_t(\gamma+\lambda)) \mathbbm{1}_{t \geq \kappa}\right).$$ Furthermore $\widehat{\eta}_\kappa$ is an infinite set containing finitely many positive points, $\eta^-_\kappa(\gamma+\lambda)\subset \widehat{\eta}_\kappa$ and $\widehat{\eta}_\kappa$ contains at least one element further right than the right edge of $\eta^-_\kappa(\gamma + \lambda)$ by the definition of $\kappa$ and the equality $\widehat{\eta}_\kappa = \eta^-_\kappa(\gamma)$. Therefore, using Lemma \ref{lemma:diffsup1} and the strong Markov property, we have 

\begin{align*}
\E \left((\widehat{r}_t - r^-_t(\gamma+\lambda)) \mathbbm{1}_{t \geq \kappa}\right) &= \P \left( t \geq \kappa \right) \E \left(\widehat{r}_t - r^-_t(\gamma+\lambda) \ | \  t \geq \kappa  \right) \\
&= \P \left( t \geq \kappa \right) \E \left(r^{\widehat{\eta}_\kappa}_{t-\kappa}(\gamma+\lambda) - r^{\eta_\kappa(\gamma+\lambda)}_{t-\kappa}(\gamma+\lambda) \ \middle| \  t \geq \kappa  \right)\\
&\geq \P \left( t \geq \kappa \right).
\end{align*}

Furthermore, if we denote by $\widehat{\kappa}$ the first time at which the rightmost element in $\left(\eta^-_t(\gamma + \lambda)\right)_{t \geq 0}$ is affected by a mark from one of the Poisson processes of rate $\lambda$, we have that $$\P \left( \kappa \leq t \right) \geq \P \left( \widehat{\kappa} \leq t \right).$$

Therefore it follows that $$\E \left( r^-_t(\gamma) \right) - \E \left( r^-_t(\gamma+\lambda) \right) \geq 1 - e^{-\lambda t}.$$

\vspace{0.4 cm}

Now for any integer $n \geq 1$
\begin{align*}
\E \left( r^-_t(\gamma) \right) - \E \left( r^-_t(\gamma+\lambda) \right) &\geq \sum_{k=1}^n \E \left[ r^-_t\left(\gamma + \frac{k-1}{n} \lambda \right)\right] - \E \left[ r^-_t\left(\gamma + \frac{k}{n} \lambda \right)\right]\\
&\geq n \left( 1 - e^{-\frac{\lambda}{n}t}\right),
\end{align*}
and the last term of these inequalities is equal to $\lambda t + o(1)$ when $n$ diverges by Taylor expansion so that by taking the limit we are left with $$ \frac{\E \left( r^-_t(\gamma) \right)}{t} - \frac{\E \left( r^-_t(\gamma+\lambda) \right)}{t} \geq \lambda.$$

We conclude by noticing that the left-hand side converges to $\alpha(\gamma) - \alpha(\gamma + \lambda)$ using the $L^1$ convergence part of Proposition \ref{prop:convtoalpha} and the assumption that $0 \leq \gamma+\lambda < \gamma_c$.
\end{proof}

\vspace{0.4 cm}

From there we obtain the following proposition, which was the purpose of this section and is mandatory to establish the results of the following section. It simply says that in the sub-critical regime the limit obtained in Proposition \ref{prop:convtoalpha} cannot be equal to $0$.

\begin{proposition} \label{prop:posalpha}
If $\gamma < \gamma_c$ then $\alpha(\gamma) > 0$.
\end{proposition}

\begin{proof}
By Lemma \ref{lemma:superlin} we have $$ \alpha(\gamma) \geq \alpha\left(\gamma + \frac{\gamma_c - \gamma}{2}\right) + \frac{\gamma_c - \gamma}{2},$$

and as $\gamma + \frac{\gamma_c - \gamma}{2} < \gamma_c$ it follows from Proposition \ref{prop:convtoalpha} that the first term in the right-hand side is greater or equal to $0$ while the second term is strictly greater than $0$, which ends the proof.
\end{proof}

\section{Sub-exponential estimates}

\label{sec:subexpest}

In this section we obtain sub-exponential estimates for the edge and the time of extinction of the dual process. Then we use these results to obtain a similar sub-exponential bound for the time correlations of the auxiliary process, on the same lines as in \cite{schonmann85}. The bounds are expressed in terms of two constants $C_1$ and $C_2$ which exact value is unimportant and will change from one result to the other. In fact the value of these constants might sometimes even change from one line to the other in the course of the same proof, in order to avoid an overload in notation.

\begin{proposition} \label{prop:exporightbord}
If $\gamma < \gamma_c$ then for any $a < \alpha$ there exists positive constants $C_1$ and $C_2$ (depending on $a$) such that for any $t \geq 0$ $$\P \left( r^-_t < a t\right) \leq C_1 e^{-C_2t}.$$ 
\end{proposition}

\begin{proof}
This result is the analogue for our system of Theorem 4 in \cite{durrett83}. The authors prove this result for another well-known interacting particle system, namely the contact process on $\Z$, using a clever construction linking the graphical characterization of the contact process (analog to the construction we gave in Section \ref{sec:themodel}) to one-dependent percolation. Nonetheless, as noticed by the authors themselves at the very end of the Section 2 of their article, this construction--and therefore the proof of their Theorem 4--can be carried out without supplementary work for a larger class of systems which includes at least nearest neighbors additive growth models. This construction is thus valid for our process as well. We let it to the reader to check that all the arguments given there works as well for our system, using previously proven results. The crucial point to carry out this construction is that there exists in the sub-critical regime some $\alpha>0$ such that $r^-_t \sim t \alpha$ as $t$ goes to $\infty$ (as proven in the previous section). The only other results needed to check the validity of their proof for our system are translation invariance, duality and the existence of spatially ergodic invariant measures with positive density to which $(\eta_t)_{t \geq 0}$ and $(\xi_t)_{t \geq 0}$ converge monotonically as $t$ goes to $\infty$ (see Section \ref{subsec:importantprop}).
\end{proof}

\vspace{0.4 cm}

In the course of the proof of Theorem \ref{th:exptinf} below we would like at some point to use the converse of the item $(iii)$ in Lemma \ref{lemma:threelemmanotdead}. Unfortunately, it suffices to consider the event in which the only dual-valid path starting at $0$ in the graphical construction immediately ends on the tips of a double arrow to see that the converse doesn't hold. A little bit of thought though reveals that this is actually the only counter-example, so that we can obtain an assertion which is close enough to the converse we need. This is the object of the next lemma. We introduce the following notation (using the jump times of Section \ref{sec:dualprocess}): $$\tilde{T} = \min \left\{\tilde{T}^*_{-1,0}, \tilde{T}^*_{0,0}, \tilde{T}^*_{1,0}, \tilde{T}^\dagger_{0,0} \right\}.$$

In other words $\tilde{T}$ corresponds to the time of the first event affecting neuron $0$ in $\left(\eta^0_t\right)_{t \geq 0}$. Moreover we define the following event $$E=\{\tilde{T} = \tilde{T^*}_{0,0}\}.$$ Now we can formulate our lemma.

\begin{lemma} \label{lemma:tau0finite}
On $E^c$, if $\sigma^0 < \infty$ then $l^+_{\sigma_0} > r^-_{\sigma_0}$.
\end{lemma}

\begin{proof}
The complementary of $E$ can be divided into the following disjoint union
\begin{equation} \label{eq:ecpart}
E^c = \left\{ \tilde{T} = \tilde{T}^\dagger_{0,0}\right\} \cup \left\{ \tilde{T} = \tilde{T}^*_{-1,0}\right\} \cup \left\{ \tilde{T} = \tilde{T}^*_{1,0}\right\}.
\end{equation}

We take care of the different cases in (\ref{eq:ecpart}) separately. First suppose that $\{\tilde{T} = \tilde{T}^\dagger_{0,0}\}$. Then obviously $\sigma^0 = \tilde{T}$. Moreover in this case $l^+_{\sigma^0} \geq 1$ and $ r^-_{\sigma^0} \leq -1$ so that $l^+_{\sigma^0} > r^-_{\sigma^0}$.

\vspace{0.4 cm}

Now let's consider the remaining cases, that is either $\tilde{T} = \tilde{T}^*_{-1,0}$ or $\tilde{T} = \tilde{T}^*_{1,0}$. By symmetry it is sufficient to treat only one them. We suppose $\tilde{T} = \tilde{T}^*_{-1,0}$. We've assumed that $\sigma^0$ is finite so that we are allowed to define $K$ to be the random variable corresponding to the index of the last active neuron before extinction. Moreover the fact that $\tilde{T} = \tilde{T}^*_{-1,0}$ implies that $\max_{0 \leq t < \sigma^0} \left| \eta^0_t \right| \geq 2$. Therefore we can define the time at which the penultimate neuron in $(\eta^0_t)_{t \geq 0}$ becomes quiescent: $$ S = \sup \{ s < \sigma^0: \left| \eta^0_s \right| = 2\}.$$ Again we shall distinguish between two separate cases: either neuron $K$ becomes quiescent because it encounters an event of the $( \tilde{T}^\dagger_{i,n})_{i\in \Z, n \in \N}$ family or because it encounters an event of the $( \tilde{T}^*_{i,n})_{i\in \Z, n \in \N}$ family.

\vspace{0.4 cm}

For the first option, we notice that by the first item of lemma \ref{lemma:threelemmanotdead} we have $r^-_s = K = l^+_s$ for any $s \in ]S,\sigma^0[$. Therefore $r^-_{\sigma^0} \leq K-1$ and $l^+_{\sigma^0} \geq K+1$, so that $l^+_{\sigma^0} > r^-_{\sigma^0}$.

\vspace{0.4 cm}

Now for the second option we let $L$ be the index of the penultimate active neuron, which becomes quiescent at time $S$. Without loss of generality we assume $L < K$ (by symmetry the case $L > K$ can be proven using the same arguments). This last case breaks down into two more sub-cases: either $L=K-1$ or $L<K-1$.

\vspace{0.4 cm}

If $L=K-1$, then neuron $L$ has $K$ as an active neighbor at time $S$, and therefore it becomes quiescent because of a `$\delta$' mark. Then there is no element of the family $(\tilde{T}^*_{L,n})_{n \in \N}$ on $]S,\sigma^0[$ as otherwise, either $K$ would not be the last particle alive, either $L$ would not be the penultimate. Therefore $K$ has no active neighbor in $(\eta^-_s)_{s \in ]S,\sigma_0[}$ and thus it becomes quiescent at time $\sigma^0$ in  $(\eta^-_s)_{s \geq 0}$ as well. Moreover by Lemma \ref{lemma:threelemmanotdead} $r^-_s = K$ for any $s \in ]S,\sigma^0[$ so that $r^-_{\sigma^0} \leq K-2$. Furthermore $l^+_{\sigma^0} \geq K$, so that we indeed have $l^+_{\sigma^0} > r^-_{\sigma^0}$ in this case. 

\vspace{0.4 cm}

Finally suppose that $L<K-1$. Define $R$ to be the time of the last event affecting $(\eta^0(t))_{t \geq 0}$ before time $S$, that is: $$R = \sup\{r<S: |\eta^0_r| \neq 2 \}.$$  Then by the second part of lemma \ref{lemma:threelemmanotdead}, for any $s \in [R,S]$, $\eta^0_s = \eta^-_s \cap [L,K]$. And as any $i \in \rrbracket L,K \llbracket$ is quiescent in $(\eta^0_s)_{s \in [R,S]}$ it shall be quiescent in $(\eta^-_s)_{s \in [R,S]}$ as well. In particular this is true for $i=K-1$. Then the conclusion follows from the same arguments as in the previous case.

\end{proof}

\vspace{0.4 cm}

Using the two previous results we are able to obtain the following precious exponential bound on the time of extinction of the dual process.

\begin{theorem} \label{th:exptinf}
If $\gamma < \gamma_c$, then there exists two positive constants $C_1$ and $C_2$ such that for any $t \geq 0$ and any finite set $A \in \mathcal{P}(\Z)$ $$\P \left(t < \sigma^A < \infty \right) \leq C_1 |A| e^{-C_2t}.$$
\end{theorem}

\begin{proof}
We first prove the result for $A = \{0\}$. Let $a$ be some real number such that $0<a<\alpha$, and let $N \geq 0$ be some integer. Then by Proposition \ref{prop:exporightbord} the following holds for some constants $C_1$ and $C_2$

\begin{equation*} \label{eq:rnan}
\P \left( r^-_n < a n \ \text{ for some } \ n \geq N\right) \leq C_1 \sum_{n \geq N} (e^{-C_2n}) = \frac{C_1}{1 - e^{- C_2}} e^{-C_2N}.
\end{equation*}

\vspace{0.4 cm}

Moreover, if $(N_t)_{t \geq 0}$ denotes an homogeneous Poisson Process of rate $1$, then an obvious coupling with $(r^-_t)_{t \geq 0}$ gives $$ \P \left(\left\{r^-_t < 0 \text{ for some } t \in [n,n+1] \right\} \cap \left\{ r^-_{n+1} > a (n+1) \right\}\right) \leq \P \left( N_1 > a(n+1)\right),$$ and by the exponential Markov inequality

\begin{equation*} \label{eq:none}
\P \left( N_1 > a(n+1)\right) \leq e^{e-1-a}e^{-an}.
\end{equation*}

\vspace{0.4 cm}

Now 
\begin{align*}
\P \left( r^-_t < 0 \text{ for some } t \geq N\right) &\leq \P \left( r^-_n < a n \ \text{ for some } \ n \geq N\right) \\
&+ \sum_{n \geq N}\P \left(\left\{r^-_t < 0 \text{ for some } t \in [m,m+1] \right\} \cap \left\{ r^-_{m+1} > a (m+1) \right\}\right).
\end{align*}

Therefore, from the inequalities above it is easy to find two constants $C_1'$ and $C_2'$ such that $$\P \left( r^-_t < 0 \text{ for some } t \geq N\right) \leq C_1' e^{-C_2'N}.$$

\vspace{0.4 cm}

And by a slight modification of the constant $C_1'$ we obtain that for any $T \in \R_+$ $$\P \left( r^-_t < 0 \text{ for some } t \geq T\right) \leq C_1' e^{-C_2'T}.$$

\vspace{0.4 cm}

By symmetry we have the same bound for $ \P \left(l^+_t > 0 \text{ for some } t \geq T\right)$, therefore 

\begin{align} \label{eq:rtlesslt}
\P \left( r^-_t < l^+_t \text{ for some } t \geq T\right) &\leq \P \left( r^-_t < 0 \text{ for some } t \geq T\right) + \P \left( l^+_t > 0 \text{ for some } t \geq T\right) \nonumber\\
&\leq 2C_1' e^{-C_2'T}.
\end{align}

\vspace{0.4 cm}

Finally
\begin{align*}
\P \left( t < \sigma^0 < \infty \right) &\leq \P \left( \left\{t < \sigma^0 < \infty\right\} \cap E^c\right) + \P \left( \{ \sigma^0>t\} \cap E\right),
\end{align*}

and the first element in the right hand side is less than $2C_1' e^{-C_2't}$ by Lemma \ref{lemma:tau0finite} and (\ref{eq:rtlesslt}), while the second element is less than $\P \left( \tilde{T}_{0,0} > t \right) = e^{-t}$, so that in the end $$\P \left( t < \sigma^0 < \infty \right) \leq 2 (C_1' + 1) e^{-\min(1,C_2')t}.$$

\vspace{0.4 cm}

Now it only remains to generalize the result to some finite initial state $A$ not necessarily equal to $\{0\}$. This is easily done using additivity, monotonicity and translation invariance
\begin{align*}
\P \left(t < \sigma^A < \infty \right) &= \P \left( \bigcup_{i \in 
A} \left\{\sigma^i > t \right\} \cap \left\{ \sigma^A < \infty \right\} \right)\\
&\leq \sum_{i \in A} \P \left( \sigma^i > t, \sigma^A < \infty\right)\\
&\leq \sum_{i \in A} \P \left( t < \sigma^i < \infty\right)\\
&= |A| \cdot \P \left(t < \sigma^0 < \infty \right)
\end{align*}

\end{proof}

\vspace{0.4 cm}
Now the only missing piece in order to prove that $(\xi_t)_{t \geq 0}$ as exponentially decaying time correlations is the following lemma, which gives a useful linear decomposition of cylinder functions in terms of indicator functions.

\begin{lemma} \label{lemma:linearcomb}
    Let $f : \mathcal{P}(\Z) \rightarrow \R$ be a cylinder function. Then there exists a family $(S_i)_{i \in \llbracket 1,n \rrbracket}$ of finite sets in $\Z$  and a family $(\lambda_i)_{i \in \llbracket 0,n \rrbracket}$ of values in $\R$ such that for any $A \in \mathcal{P}(\Z)$ $$ f(A) = \lambda_0 + \sum_{i=1}^n \lambda_i \mathbbm{1}_{A \cap S_i \neq \emptyset}.$$
\end{lemma}

\begin{proof}
Let $B \subset \Z$ be the support of $f$, and let $(B_i)_{i \in \llbracket 1, n\rrbracket}$ be the family of all the subsets of $B$. Then, for any $A \in \mathcal{P}(\Z)$
\begin{equation} \label{eq:initialinear}
    f(A) = f(A \cap B) = \sum_{i=1}^n f(B_i) \mathbbm{1}_{A \cap B = B_i}.
\end{equation}

Moreover, for any $A \in \mathcal{P}(\Z)$ the following identity holds \begin{equation} \label{eq:recformula}
    \mathbbm{1}_{A \cap B = B_i} = \mathbbm{1}_{A \cap B \subset B_i} - \sum_{\substack{ B_j \subset B_i \\ B_j \neq B_i}} \mathbbm{1}_{A \cap B=B_j},
\end{equation} and if $B_j = \emptyset$  then
\begin{equation} \label{eq:singleton}
    \mathbbm{1}_{A \cap B = B_j} = \mathbbm{1}_{A \cap B \subset B_j}.
\end{equation}

\vspace{0.4 cm}

Therefore, applying (\ref{eq:recformula}) recursively to itself---that is, to the elements in the sum of the right hand side---and using (\ref{eq:singleton}) to conclude the recursion, it follows from (\ref{eq:initialinear}) that $f$ can be expressed as a finite linear combination of indicators of the form $\mathbbm{1}_{\text{ }\kern-0.10em \bigcdot \text{ } \cap B \subset B_i}$ as well. Then for any $i \in \llbracket 1,n\rrbracket$ we define $S_i = B_i^c \cap B$ and the conclusion follows from the fact that for any $A \in \mathcal{P}(\Z)$ the following holds $$\mathbbm{1}_{A \cap B \subset B_i} = 1 - \mathbbm{1}_{A \cap S_i \neq \emptyset}.$$

\end{proof}

\begin{theorem} \label{th:timecorrelation}
Let $f: \mathcal{P}(\Z) \rightarrow \R$ be a cylinder function. If $\gamma < \gamma_c$ then there exists positive constants $C_1$ and $C_2$ (depending on $f$) such that for any $s,t \in \R^+$ $$|\Cov \left( f(\xi_t), f(\xi_s)\right)| \leq C_1 e^{-C_2 |t-s|}.$$
\end{theorem}

\begin{proof}
Let $(S_i)_{i \in \llbracket 1,n \rrbracket}$ and $(\lambda_i)_{i \in \llbracket 1,n \rrbracket}$ be as in Lemma \ref{lemma:linearcomb} so that for any $t \geq 0$ $$ f(\xi_t) = \sum_{i=1}^n \lambda_i \mathbbm{1}_{\xi_t \cap S_i \neq \emptyset}.$$

Therefore, for any $s \geq 0$ and $t \geq 0$ 

$$ \left| \Cov \left(f(\xi_t), f(\xi_s)\right) \right| \leq \sum_{i=1}^n \sum_{j=1}^n \left| \lambda_i \lambda_j \right| \cdot \left| \Cov \left(\mathbbm{1}_{\xi_t \cap S_i \neq \emptyset},\mathbbm{1}_{\xi_s \cap S_j \neq \emptyset} \right)\right|.$$

Thus it is sufficient to prove that for any finite sets $A,B \in \mathcal{P}(\Z)$ there exists some positive constants $C_1$ and $C_2$ such that for any $s \geq 0$ and $t \geq 0$ 

$$|\Cov \left( \mathbbm{1}_{\xi_t \cap A \neq \emptyset}, \mathbbm{1}_{\xi_s \cap B \neq \emptyset}\right)| \leq C_1 e^{-C_2 |t-s|}.$$

\vspace{0.4 cm}

Without loss of generality we assume $t \leq s$. Let $F$ be the event that there exists a valid path from $\Z \times t$ to $B \times s$. Moreover let $G = \{ \xi_t \cap A \neq \emptyset \}$ and $H = \{ \xi_s \cap B \neq \emptyset \}$. We have

$$\left|\Cov \left( \mathbbm{1}_{\xi_t \cap A \neq \emptyset}, \mathbbm{1}_{\xi_s \cap B \neq \emptyset}\right)\right| = \left| \P \left(G \cap H\right) - \P \left(G\right) \P \left( H\right)\right|.$$

\vspace{0.4 cm}

Furthermore it is clear from the graphical construction that $H \subset F$ so that we can replace $H$ by $H \cap F$ in the equation above. Also events depending on disjoint regions of the graph $\mathcal{G}$ are independent so that $P\left( G \cap F\right) = P\left( G \right) \P\left(F\right)$. Therefore 

\begin{align*}
\left| \Cov \left( \mathbbm{1}_{\xi_t \cap A \neq \emptyset}, \mathbbm{1}_{\xi_s \cap B \neq \emptyset}\right) \right| &= \left| \P \left(G \cap H \cap F\right) - \P(G \cap F) + \P(G)\P(F) - \P \left(G\right) \P \left( H \cap F\right)\right|\\
&=\left|  \left[\P(G \cap F) - \P \left(G \cap H \cap F\right) \right] - \P(G) \left[ \P(F) - \P(H \cap F)\right]\right|\\
&= \left| \P \left(G \cap F \cap H^c\right) - \P(G)  \P(F \cap H^c)\right|\\
&= \P\left(F \cap H^c \right) \cdot \left| \P(G) - \P \left(G | F \cap H^c\right) \right|\\
&\leq \P\left(F \cap H^c \right).
\end{align*}

\vspace{0.4 cm}

And by the backward version of the dual process, 

$$\P\left(F \cap H^c\right) = \P \left(s- t < \sigma^B < s\right) \leq \P \left(s - t < \sigma^B < \infty\right).$$

\vspace{0.4 cm}

Then use Theorem \ref{th:exptinf} to conclude.

\end{proof}

\vspace{0.4 cm}

\section{Proof of the main theorem}

\label{sec:mainthm}

Before proving Theorem \ref{thm:main} per se we shall prove a similar result concerning the auxiliary process, which is the object of Section \ref{subsec:timeavaux}.
Before doing so we introduce the finite  versions of the auxiliary process in Section \ref{sec:finsemifinaux}, as well as the semi-infinite versions, and explain how, together with the infinite system, they all relate to each other. We then state important results about the time of extinction of the finite version in Section \ref{sec:resulttimeext}. In Section \ref{sec:mainproof} we then proceed to the main proof.

\vspace{0.4 cm}

\subsection{Finite and semi-infinite auxiliary processes}

\label{sec:finsemifinaux}

The \textbf{finite auxiliary processes} are the auxiliary processes of the finite versions of the system, with which Theorem \ref{thm:main} is concerned, that is, the cases $I_n = \llbracket -n, n \rrbracket$ and $\V_{n,i} = \{-i,i\} \cap \llbracket -n, n \rrbracket$. These are the processes with state space $\{0,1\}^{2n + 1}$ and which dynamic is given by generator (\ref{generator}) when you replace $\Z$ by $I_n$ under the summation, and where the maps (\ref{spikesmap}) are defined with respect to $\V_{n,i}$ instead of $\{i-1,i+1\}$. We write $\big(\xi_n(t)\big)_{t \geq 0}$ for these finite versions. Moreover we define the extinction time of the finite process, which is always finite by standard considerations about irreducible Markov processes on a finite state space. For any $A \in \mathcal{P}(\llbracket -n,n \rrbracket)$ let $$\tau_n^A = \inf \{t \geq 0: \xi_n^A(t) = \emptyset\}.$$

\vspace{0.4 cm}

We also introduce the \textbf{semi-infinite auxiliary processes}, corresponding to the case in which the neurons are indexed on $I_{[n,\infty[} = \integerninfty$ or $I_{]-\infty,n]} = \integerinftyn$ for some $n \in \Z$, the set of presynaptic neurons being given as usual by $\{i-1,i+1\} \cap I_{[n,\infty[}$ or $\{i-1,i+1\} \cap I_{]-\infty,n]}$. These processes are denoted $\big(\xi_{[n,\infty[}(t)\big)_{t \geq 0}$ and $\big(\xi_{]-\infty,n]}(t)\big)_{t \geq 0}$ respectively. These semi-infinite processes have the same general properties as their infinite counterpart; in particular there exists an upper-invariant measure for each of them, giving mass $0$ to the empty-set in the sub-critical regime (see \cite{andre19} for details). We write $\mu_{[0,\infty[}$ for the upper-invariant measure of $\big(\xi_{[0,\infty[}(t)\big)_{t \geq 0}$, which is the only one that will be needed in the proofs below. 

\vspace{0.4 cm}

A graphical construction similar to the one introduced in Section \ref{subsec:graphconsaux} can obviously be given for the finite and semi-infinite processes, using only the Poisson processes of one of the sub-diagrams $\llbracket -n, n \rrbracket \times \R_+$, $\integerninfty \times \R_+$ or $\integerinftyn \times \R_+$. This allows a simple and very useful coupling between the finite, semi-infinite and infinite processes, which gives us the following lemma.

 \begin{lemma} \label{lemma:finiteandinfinite}
	For any $t \geq 0$ the following holds on $\{\tau_n > t\}$
	
	\begin{enumerate}[(i)]
		\item $\xi_n(t) = \xi(t) \cap [\min \xi_n(t), \max \xi_n(t)]$,
		\item $\min \xi_n(t) = \min \xi_{[-n,\infty[}(t)$ and $\max \xi_n(t) = \min \xi_{]-\infty,n]}(t)$.
	\end{enumerate}
	
\end{lemma}

\vspace{0.4 cm}

The proof, which is pretty straightforward using the natural coupling permitted by the graphical construction, can be found in \cite{andre19}.

\vspace{0.4 cm}

\subsection{Results about the time of extinction of the finite process}

\label{sec:resulttimeext}

Now let us state a fundamental result about the time of extinction of $(\xi_n(t))_{t \geq 0}$ which was proven in \cite{andre19}, and which is the first of the two properties characterizing metastable systems we've evoked in the introduction. Notice that while this result is stated in terms of the extinction time of the finite auxiliary process--that is the time at which the last neuron becomes quiescent--it trivially implies the same result for the time of the last spike, even if these two random times are not formally equals. This result will be used repeatedly from now on up to the end of this article.

\begin{theorem} \label{th:meta1}
	Suppose $0 < \gamma < \gamma_c$. Then we have the following convergence
	
	$$\frac{\tau_n}{\E (\tau_n)} \overset{\mathcal{D}}{\underset{n \rightarrow \infty}{\longrightarrow}} \mathcal{E} (1),$$ where $\mathcal{D}$ denotes a convergence in distribution.
\end{theorem}

\begin{remark}
	In \cite{andre19} this result was proven only for $\gamma < \gamma'_c$, where $\gamma'_c$ is the critical value for the semi-infinite process. Nonetheless, with the construction invoked in the proof of Proposition \ref{prop:exporightbord} the parameter $p$ of the percolation process can be taken as close to $1$ as needed, while the coupling with $(\xi_t)_{t \geq 0}$ remains valid even if the value of $\gamma$ is kept fixed. But if $p$ is big enough then with positive probability there is a path in the percolation structure which never goes further left than $0$ on the horizontal abscissa. An immediate consequence is that $\gamma_c = \gamma'_c$, so that we don't need the caveat anymore.
\end{remark}

\vspace{0.4 cm}

The second result tells us that in the sub-critical regime the expectation of the time of extinction grows strictly faster than $n$ itself as $n$ grows. While this result will mostly be needed in the last section of this paper, for the moment it gives us a rigorous proof of a fact that might seems anyway evident: it ensures us that $\E (\tau_n) \rightarrow \infty$ as $n$ grows, which in turns guarantees that there is no problem in choosing a $(R_n)_{n \geq 0}$ satisfying the conditions of Theorem \ref{thm:main}.

\vspace{0.4 cm}

\begin{proposition} \label{thm:expectsuperlinear}
	Suppose $0 < \gamma < \gamma_c$. Then $n/\E(\tau_n) \rightarrow 0$ as $n$ goes to $\infty$.
\end{proposition}

\begin{proof}
	In the following we write $\frac{n}{2}$ while we should sometimes write $\left\lfloor \frac{n}{2}\right\rfloor$ in order to avoid a useless overload in the notation. By duality $$ \P \left( \xi^{\llbracket 0,\frac{n}{2} \rrbracket}_t \neq \emptyset \right) \underset{t \rightarrow \infty}{\longrightarrow} \widetilde{\mu} \left( A : A \cap \left[0,\frac{n}{2}\right] \neq \emptyset \right).$$
	
	Therefore, if we write $E_n \mydef \left\{ \xi^{\llbracket 0,\frac{n}{2} \rrbracket}_t \neq \emptyset \text{ for all } t \geq 0 \right\}$ then $E_n = \bigcap_{t \geq 0} \left\{\xi^{\llbracket 0,\frac{n}{2} \rrbracket}_t \neq \emptyset\right\}$ and thus $$ \P \left( E_n \right) = \widetilde{\mu} \left( A : A \cap \left[0,\frac{n}{2}\right] \neq \emptyset \right) \underset{n \rightarrow \infty}{\longrightarrow} 1.$$
	
	Now we define $$\kappa_n \mydef \inf \left\{t > 0: \{-n,n\} \cap \xi^{\llbracket0,\frac{n}{2}\rrbracket}_t \neq \emptyset\right\}.$$
	
	From Lemma \ref{lemma:threelemmanotdead} (and translation invariance) we have $r^{\llbracket 0, \frac{n}{2}\rrbracket}_t = r^{\rrbracket - \infty, \frac{n}{2}\rrbracket}_t$ on $E_n$, and therefore by Proposition \ref{prop:convtoalpha} and Proposition \ref{prop:posalpha} it follows that $r^{\llbracket 0, \frac{n}{2}\rrbracket}_t \rightarrow \infty$ almost surely on $E_n$ when $t$ diverges. As a consequence, on $E_n$, the stopping time $\kappa_n$ is almost surely finite and $\tau_n \geq \kappa_n$. Moreover, if $\left( W_k\right)_{k \geq 1}$ is a family of independent and identically distributed random variables such that $W_1 \leadsto \mathcal{E}(1)$, then an obvious coupling gives that $$\kappa_n \geq \sum_{k=1}^{n/2} W_k.$$
	
	Therefore $$ \P \left( \tau_n < \frac{n}{4}, \ E_n\right) \leq \P \left( \frac{1}{n/2} \sum_{k=1}^{n/2} W_k > \frac{1}{2}\right).$$

	The right hand side in the inequality above goes to $0$ as $n$ diverges by the law of large numbers. Thus, as $E_n$ has probability one asymptotically it implies that $\P \left( \tau_n < n/4\right)$ goes to $0$ when $n$ goes to $\infty$. Then dividing by $\E(\tau_n)$ in both side of the inequality and using Theorem \ref{th:meta1} we obtain the result.
	
\end{proof}

\vspace{0.4 cm}

\subsection{Convergence of the times averages for the auxiliary process}
\label{subsec:timeavaux}

For any $t,R \in \R_+$ and $n \geq 0$ we define the following measure on $\mathcal{P}(\Z)$ $$A_R^n(t, \ \bigcdot \ ) = \frac{1}{R} \int_t^{t+R} \mathbbm{1}_{\xi_n(s) \in \ \bigcdot \ } ds.$$

\vspace{0.4 cm}

As usual, for any measurable function $f : \mathcal{P}(\Z) \rightarrow \R$ we write $A_R^n(t,f)$ for the integral of $f$ with respect to $A_R^n(t, \ \bigcdot \ )$. It is straightforward from basic measure-theoretic considerations to see that actually $$ A_R^n(t,f) = \frac{1}{R} \int_t^{t+R} f(\xi_n(s)) ds.$$

\vspace{0.4 cm}

Then the following result holds.

\begin{theorem} \label{thm:mainaux}

Let $(R_n)_{n \geq 0}$ be an increasing sequence of positive real numbers satisfying the same conditions as in Theorem \ref{thm:main}. Then, for any cylinder function $f : \mathcal{P} (\Z) \mapsto \R$ we have, for any $t \geq 0$

$$ A_{R_n}^n \left( t,f\right) \overset{\P}{\underset{n \rightarrow \infty}{\longrightarrow}} \mu(f).$$
\end{theorem}

\vspace{0.4 cm}

The two main ingredients are the following lemmas, which we prove now, before entering the proof of Theorem \ref{thm:mainaux}. These are essentially a consequence of the weak convergence of $(\xi_t)_{t \geq 0}$, together with Chebyshev's inequality and the exponential decay of the time correlations.

\vspace{0.4 cm}

\begin{lemma} \label{lemma:mainaux1}
	Let $\epsilon > 0$ and $f : \mathcal{P}(\Z) \rightarrow \R$ be a cylinder function. Then there exists a constant $C>0$ (which depends only on $f$ and $\epsilon$) such that for any fixed $t \geq 0$ and for big enough $R \in \R_+$ the following holds $$\P \left(\left| \frac{1}{R} \int_{t}^{t + R} f(\xi_s) ds - \mu(f)\right| > \epsilon\right) \leq \frac{C}{R}.$$
\end{lemma}

\begin{proof}
	Fix some $t \geq 0$. As $\mu$ is defined as the weak limit of $(\xi_t)_{t \geq 0}$ there exists some $t_0$ such that for any $t \geq t_0$ $$\left| \E \left( f(\xi_t) \right) - \mu(f) \right| < \frac{\epsilon}{2}.$$
	
	Then, if $t \geq t_0$ we have
	\begin{equation} \label{eq:diffespmu}
		\left|\frac{1}{R} \int_{t}^{t + R} \E \left( f(\xi_s) \right)ds - \mu(f) \right| < \frac{\epsilon}{2}. 
	\end{equation}
	
	\vspace{0.4 cm}
	
	If $t<t_0$ then (assuming $t + R > t_0$) $$  \left|\frac{1}{R} \int_{t}^{t+R} \E \left( f(\xi_s) \right)ds - \mu(f) \right| \leq \frac{1}{R} \int_{t}^{t_0} \left| \E \left( f(\xi_s) \right) - \mu(f)\right| ds + \frac{1}{R} \int_{t_0}^{t + R} \left| \E \left( f(\xi_s) \right) - \mu(f)\right| ds,$$ and as the first element in the sum of the right hand side goes to $0$ as $R$ grows while the second element is strictly less than $\epsilon / 4$, the following must holds for big enough $R$
	\begin{equation} \label{eq:boundxtr}
		\P \left( \left| \frac{1}{R} \int_{t}^{t + R} f(\xi_s) ds - \mu(f)\right| > \epsilon \right) \leq \P \left( \left| \frac{1}{R} \int_{t}^{t + R} f(\xi_s) ds - \frac{1}{R} \int_{t}^{t + R}\E \left(f(\xi_s)\right) ds \right| > \frac{\epsilon}{2}\right).
	\end{equation}
	
	\vspace{0.4 cm}
	
	Writing $ X_{t,R} \mydef \frac{1}{R} \int_{t}^{t + R} f(\xi_s)ds$ and using Fubini's theorem, (\ref{eq:boundxtr}) becomes
	
	\begin{equation} 
		\P \left( \left| \frac{1}{R} \int_{t}^{t + R} f(\xi_s) ds - \mu(f)\right| > \epsilon \right) \leq \P \left( \left| X_{t,r} - \E \left(X_{t,R}\right) \right| > \frac{\epsilon}{2}\right).
	\end{equation}
	
	\vspace{0.4 cm}
	
	Again, Fubini's theorem gives $$\E\left(X_{t,R}^2\right) = \frac{1}{R^2} \int_{t}^{t+R} \int_{t}^{t + R} \E \left( f(\xi_u)f(\xi_v)\right) du dv,$$ and $$\E\left(X_{t,R}\right)^2 = \frac{1}{R^2} \int_{t}^{t+R} \int_{t}^{t+R} \E \left( f(\xi_u)\right) \E \left(f(\xi_v)\right) du dv,$$
	and therefore from Theorem \ref{th:timecorrelation} we obtain the following
	\begin{align*}
		\Var(X_{t,R}) &= \frac{1}{R^2} \int_{t}^{t+R} \int_{t}^{t+R} \Cov \left( f(\xi_u), f(\xi_v)\right) du dv\\
		&\leq \frac{1}{R^2} \int_{t}^{t+R} \int_{t}^{t+R} C_1 e^{-C_2 |u - v|} du dv.
	\end{align*}
	
	\vspace{0.4 cm}
	
	Then by a simple change of variable, for any $u \in [t, t + R]$ \begin{align*}
		\int_{t}^{t + R} C_1 e^{-C_2 |u - v|} dv &= C_1 \left[ \int_0^{u - t} e^{-C_2 v} dv + \int_0^{t+R-u} e^{-C_2 v} dv \right] \\
		&\leq 2 C_1 \int_0^{\infty} e^{-C_2 u} du.
	\end{align*}
	
	Thus $$ \Var \left( X_{t,R} \right) \leq \frac{2 C_1}{C_2 R}.$$ 
	
	\vspace{0.4 cm}
	
	Then Chebyshev's inequality and (\ref{eq:boundxtr}) imply that $$ 	\P \left( \left| \frac{1}{R} \int_{t}^{t + R} f(\xi_s) ds - \mu(f)\right| > \epsilon \right) \leq\frac{8 C_1}{C_2 \epsilon^2 R}.$$
	
\end{proof}

\vspace{0.4 cm}

\begin{lemma} \label{lemma:mainaux2}
	Let $\epsilon > 0$ and $f : \mathcal{P}(\Z) \rightarrow \R$ be a cylinder function. Then there exists a positive constant $C$ (which depends only on $f$ and $\epsilon$) such that if $l$ is big enough, then for any fixed $t \geq 0$ and for any $R \in \R^+$ the following holds $$\P \left[\frac{1}{R}\int_{t}^{t + R} \mathbbm{1}_{\xi_{[0,\infty[}(s)  \cap  [0,l] = \emptyset} \ ds > \epsilon\right] \leq \frac{C}{R}.$$
\end{lemma}

\begin{proof}
	This lemma is easily obtained by similar arguments as in the proof of Lemma \ref{lemma:mainaux1}. For any $l \geq 0$, the function $h_l : \xi \mapsto \mathbbm{1}_{\xi  \cap  [0,l] = \emptyset}$ is a decreasing function, so that by stochastic monotonicity: $$\frac{1}{R} \int_{t}^{t+R} \E \left( h_l(\xi_{[0,\infty[}(s))\right) ds \leq \mu_{[0,\infty[} \left(h_l\right).$$
	
	Moreover the fact that $\mu_{[0,\infty[}$ gives mass $0$ to the empty set in the sub-critical regime implies that $$\mu_{[0,\infty[}(h_l) = \mu_{[0,\infty[} \left( \{A: A\cap[0,l] = \emptyset\} \right) \underset{l \rightarrow \infty}{\longrightarrow} 0.$$
	
	\vspace{0.4 cm}
	
	Therefore, for big enough $l$ we have $$\P \left[\frac{1}{R}\int_{t}^{t + R} h_l(\xi_{[0,\infty[}(s)) \ ds > \epsilon \right] \leq \P \left[\frac{1}{R} \int_{t}^{t+R} h_l(\xi_{[0,\infty[}(s)) - \E \left( h_l(\xi_{[0,\infty[}(s))\right) ds > \epsilon\right].$$
	
	\vspace{0.4 cm}
	
	Then the conclusion comes from the same arguments\footnote{Notice that to apply the same arguments as in the proof of Lemma \ref{lemma:mainaux1} we need to have the exponential decay of the time correlations for the semi-infinite processes, while formally Theorem \ref{th:timecorrelation} is proven only for the infinite process. Nonetheless one can prove an analogous result for the semi-infinite processes by the same arguments as those of Theorem \ref{th:timecorrelation}.} as in the proof of Lemma \ref{lemma:mainaux1} (from Eq. (\ref{eq:boundxtr}) to the end), using the fact that $h_l$ is a cylinder function.
\end{proof}

\vspace{0.4 cm}

\begin{proof}[Proof of Theorem \ref{thm:mainaux}]
	
	Fix some $t \geq 0$. First notice that using Theorem \ref{th:meta1} and the hypothesis on $(R_n)_{n \geq 0}$ (as well as the fact that $\E(\tau_n) \rightarrow \infty$) the following holds $$ \P \Big( t + R_n < \tau_n \Big) \underset{n \rightarrow \infty}{\longrightarrow} 1.$$
	
	Hence, writing $\Omega_n = \{t + R_n < \tau_n\}$, it will be enough to prove that, for any $\epsilon > 0$ $$ \P \left( \left| A_{R_n}^n(t,f) -  \mu(f) \right| > \epsilon, \ \Omega_n \right) \underset{n \rightarrow \infty}{\longrightarrow} 0.$$
	
	We have 
	
	\begin{align*} \P \Big( \left| A_{R_n}^n(t,f) -  \mu(f) \right| > \epsilon, \ \Omega_n \Big) &\leq \P \left( \left|  A_{R_n}^n(t,f) - \frac{1}{R_n} \int_t^{t + R_n} f(\xi_s) ds \right| > \frac{\epsilon}{2}, \ \Omega_n \right)\\
	&+ \P \left( \left| \frac{1}{R_n} \int_t^{t + R_n} f(\xi_s) ds - \mu(f) \right| > \frac{\epsilon}{2} \right)
	\end{align*}

	\vspace{0.4 cm}

	Moreover, for any $l \in \N$ and for $n$ big enough, using the fact that the support of $f$ has to lie in $[-n+l,n-l]$ when $n$ is big, Lemma \ref{lemma:finiteandinfinite} yields 
	
	\begin{equation} \label{eq:bordersupport}
		\{\min \xi_n(t) < -n + l, \max \xi_n(t) > n-l \} \subset \{f(\xi_n(t)) = f(\xi_t)\}.
	\end{equation}

	\vspace{0.4 cm}
	
	For any $l \in \N$ we let $h^n_l: \mathcal{P}(\Z) \rightarrow \R$ and $g^n_l : \mathcal{P}(\Z) \rightarrow \R$ be defined for any $\xi \in \mathcal{P}(\Z)$ by $$h^n_l(\xi) = \mathbbm{1}_{\xi  \cap  [-n, -n+l] = \emptyset},$$ and $$g^n_l(\xi) = \mathbbm{1}_{\xi  \cap  [n-l, n] = \emptyset}.$$
	
	\vspace{0.4 cm}
	
	Now, using (\ref{eq:bordersupport}) and the fact that $|f(\xi_n(t)) - f(\xi(t))| < 2\|f\|_\infty$ we have
	
	\begin{equation} \label{eq:lesstwonorms}
		\left|A_{R_n}^N (t,f) - \frac{1}{R_n} \int_{t}^{t + R_n} f(\xi_s)ds \right| \leq 2\|f\|_\infty \frac{1}{R_n} \int_{t}^{t+R_n} h^n_l(\xi_n(s)) + g^n_l(\xi_n(s))ds 
	\end{equation}
	
	\vspace{0.4 cm}
	
	Furthermore, for any $l<n$,  and for any $0 \leq s < \tau_n$ Lemma \ref{lemma:finiteandinfinite} yields:
	
	\begin{equation}
		h^n_l\left(\xi_n(s)\right) = h^n_l \left( \xi_{[-n,+\infty[}(s)\right) \\
		\text{ and }\\
		g^n_l\left(\xi_n(s)\right) = g^n_l \left( \xi_{]-\infty,n]}(s)\right).
	\end{equation}
	
	\vspace{0.4 cm}
	
	Therefore on $\Omega_n$ the right hand side in (\ref{eq:lesstwonorms}) is less than $$2\|f\|_\infty \int_{t}^{t + R_n} h^n_l(\xi_{[-n,+\infty[}(s)) + g^n_l(\xi_{]-\infty,n]}(s))ds.$$
	
	\vspace{0.4 cm}
	
	Then, for any $l \in \N$ and for $n$ big enough the following holds (using translation invariance and symmetry)
	
	\begin{align*}
		\P \Big( \left| A_{R_n}^n(t,f) -  \mu(f) \right| > \epsilon, \ \Omega_n \Big) &\leq 2 \P \left[ \frac{1}{R}\int_{t}^{t + R} \mathbbm{1}_{\xi_{[0,\infty[}(s)  \cap  [0,l] = \emptyset} \ ds > \frac{\epsilon}{8 \|f\|_\infty} \right]\\
		&+ \P \left[ \left| \frac{1}{R} \int_{t}^{t + R} f(\xi_s) ds - \mu(f)\right| > \frac{\epsilon}{2} \right].
	\end{align*}

	\vspace{0.4 cm}
	
	Then the result follows from Lemma \ref{lemma:mainaux1} and Lemma \ref{lemma:mainaux2}.
\end{proof}

\vspace{0.4 cm}
	
\subsection{Proof of the main theorem}

\label{sec:mainproof}

We can now turn to the proof of the main Theorem. An important element of this proof is the following lemma, which is the last result we prove before engaging into the proof of this theorem. It is a coarser version of Theorem \ref{th:timecorrelation}, but for the finite version of our system, the exponential decay being tempered by two other terms which converge to zero when $n$ grows. While we conjecture that the exponential decay actually holds without these two terms, we only prove this weaker version as it is sufficient for our needs.  

\begin{lemma} \label{lemma:covfinite}
	Let $s,r \in \R_+$ and let $i \in \Z$. Then there exists two positive constants $C_1$ and $C_2$ such that for any $n \geq 0$ $$\left| \Cov \left( \mathbbm{1}_{\xi_n(s) \cap \{i\} \neq \emptyset}, \mathbbm{1}_{\xi_n(t) \cap \{i\} \neq \emptyset}\right) \right| \leq C_1 e^{-C_2 |t-s|} + \P \left( \tau_n < \max (s,t) \right) + \epsilon_n,$$
	
	where $\epsilon_n$ is some positive quantity satisfying $\epsilon_n \underset{n \rightarrow \infty}{\longrightarrow} 0.$
\end{lemma}

\begin{proof}
	Without loss of generality we assume $t \leq s$. Moreover we let $ G_n = \{\xi_n(t) \cap \{i\} \neq \emptyset\}$ and $ H_n = \{\xi_n(s) \cap \{i\} \neq \emptyset\}$. We also define the event $F$ that there is a valid path on the graphical construction (of the infinite process) from $\Z \times t$ to $\{i\} \times s$. Then by the same arguments as in the proof of Theorem \ref{th:timecorrelation} we obtain $$ \left| \Cov \left( \mathbbm{1}_{\xi_n(s) \cap \{i\} \neq \emptyset}, \mathbbm{1}_{\xi_n(t) \cap \{i\} \neq \emptyset}\right) \right| \leq \P \left( F \cap H_n^c \right).$$
	
	Then $$\P \left( F \cap H_n^c \right) \leq \P \left( s-t < \sigma_i < s\right) + \P \left(\tau_n \leq s \right) + \P \left( \{ \sigma_i \geq s \} \cap H_n^c \cap \{\tau_n>s\}\right).$$
	
	As in the proof of \ref{th:timecorrelation} $\P \left( s-t < \sigma_i < s\right) \leq \P \left( s-t < \sigma_i < \infty\right) \leq C_1 e^{-C_2 |t-s|}$, so it only remains to show that the second term in the sum above goes to $0$ as $n$ diverges.  
	
	\vspace{0.4 cm}
	
	If $\sigma_i \geq s$ then $ \{j\} \times 0 \longrightarrow \{i\} \times s$ for some $j \in \Z$. But, if $j \in \{-n, \ldots, n \}$, then on $H_n^c$ this path has to cross at least one of the two frontiers $\{n\} \times \R$ and $\{-n\} \times \R$ at some point, as otherwise there would be a valid path from $\{-n, \ldots n\}$ to $(i,s)$ that never escapes $\{-n, \ldots n\} \times \R^+$, which is not allowed on $H_n^c$. If $j \notin \{-n, \ldots n\}$ of course the path crosses one of the two frontiers as well. Consequently, in both cases we can define $s' \leq s$ to be the last time of crossing, and we have either $(n,s') \longrightarrow (i,s)$ or $(-n,s') \longrightarrow (i,s)$. But in the first case $\xi_n(s) \cap \{i, \ldots n\} = \emptyset$, as otherwise there would be a valid path crossing the path from $(s',n)$ to $(s,i)$ and the concatenation would then be a valid path from $\{-n, \ldots n \} \times {0}$ to $(i,s)$ never escaping $\{-n, \ldots n\} \times \R^+$ (which again is not allowed on $H_n^c$). Similarly we have $\xi_n(s) \cap \{-n, \ldots i\} = \emptyset$ in the other case.
	
	\vspace{0.4 cm}
	
	From the discussion above it follows that 
	\begin{align*}
	\P \left( \{ \sigma_i \geq s \} \cap H_n^c \cap \{\tau_n>s\}\right) &\leq \P \left( \xi_n(s) \cap \{i, \ldots n \} = \emptyset, \tau_n > s\right)\\
	&+ \P \left( \xi_n(s) \cap \{-n, \ldots i \} = \emptyset, \tau_n > s\right).
	\end{align*}
	
	\vspace{0.4 cm}
	
	Then, using Lemma \ref{lemma:finiteandinfinite} this becomes 
	\begin{align*}
	P \left( \{ \sigma_i \geq s \} \cap H_n^c \cap \{\tau_n>s\}\right) &\leq \P \left( \xi_{[0, +\infty]} (s) \cap [0,n-i] = \emptyset \right) + \P \left(\xi_{[0, +\infty]} (s) \cap [0,n+i] = \emptyset \right)\\
	&\leq 2 \P \left( \xi_{[0, +\infty]} (s) \cap [0,\min(n-i,n+i)] = \emptyset \right),
	\end{align*}

	and by stochastic monotonicity this is less than $$\epsilon_n \mydef 2 \mu_{[0,+\infty[} \left( A: A \cap [0, \min (n-i,n+i)] = \emptyset \right).$$

	\vspace{0.4 cm}
	
	Finally $\epsilon_n \underset{n \rightarrow \infty}{\longrightarrow} 0$ as $\mu_{[0,+\infty[}$ gives mass $0$ to $\emptyset$.	
\end{proof}

\vspace{0.4 cm}

With this preliminary completed, we can turn to the proof of the main theorem. 

\vspace{0.4 cm}

\begin{proof}[Proof of Theorem \ref{thm:main}]
We fix $t$ and $F$ and we let $\epsilon > 0$. We aim to show that $$ \P \left( \left| \widehat{N}^n_{R_n} \left(t,F\right) - |F| \cdot \rho \ \right| > \epsilon \right) \underset{n \rightarrow \infty}{\longrightarrow} 0.$$

\vspace{0.4 cm}

For any $\xi \in \{0,1\}^\Z$ we write $S_F (\xi) = \sum_{i \in F} \xi_i$. Then

\begin{align*}\P \left( \left| \widehat{N}^n_{R_n} \left(t,F\right) -  |F| \cdot \rho \ \right| > \epsilon \right) &\leq \P \left( \left| \widehat{N}^n_{R_n} \left(t,F\right) - A^n_{R_n} \left(t,S_F\right) \right|  > \frac{\epsilon}{2}\right)\\
&+ \P \left( \Big| A^n_{R_n} \left(t,S_F\right) - |F| \cdot \rho \ \Big| > \frac{\epsilon}{2} \right)
\end{align*}

\vspace{0.4 cm}

We already know that the second element in the sum above goes to $0$ by Theorem \ref{thm:mainaux}, so that it only remains to show that the first element goes to $0$ as well. To do so we define, for any $i \in F$ and $n  \geq 0$, the following random set on $\R^+$ $$I_{n,i} \mydef \{s: \ t\leq s\leq t+R_n  \text{ and } \xi_{n,i}(s) = 1\}.$$ 

\vspace{0.4 cm}

Then the following holds 
\begin{align*}
 \left| \widehat{N}^n_{R_n} \left(t,F\right) - A^n_{R_n} \left(t,S_F\right) \right| &\leq \frac{1}{R_n} \sum_{i \in F} \left| N_i([t,t+R_n]) - \int_t^{t+R_n} \xi_{n,i}(s) ds\right|\\
 &= \frac{1}{R_n} \sum_{i \in F} \left| N_i(I_{n,i}) - \lambda (I_{n,i})\right|,
\end{align*}

where $\lambda(I_{n,i})$ denotes the Lebesgue's measure of the set $I_{n,i}$ (there is no problem of measurability as $I_{n,i}$ is almost surely a finite union of intervals). Above we've used the fact that there can not be any spike on the set $[t, t+R_n] \cap {I_{n,i}}^c$. Furthermore, for any fixed $I_{n,i}$, $N_i(I_{n,i})$ corresponds to the numbers of atoms of a Poisson process of intensity $1$ on a Borel set of length $\lambda(I_{n,i})$, so that we immediately obtain

\begin{align} \label{eq:espcondeqlambda}
	\E \left[N_i(I_{n,i}) | \lambda(I_{n,i}) \right	] = \lambda(I_{n,i})
\end{align}

\vspace{0.4 cm}

It follows that

\begin{align*}
	\P \left[ \frac{1}{R_n} \sum_{i \in F} \left| N_i(I_{n,i}) - \lambda (I_{n,i})\right| > \epsilon \right]	&\leq \sum_{i \in F} \P \left[ \big| N_i(I_{n,i}) - \E\left[N_i(I_{n,i})\big] \right| > \frac{\epsilon R_n}{2 |F|} \right]\\
	&+ \sum_{i \in F} \P \left[ \big| \E\left[N_i(I_{n,i})\right] - \E\left[N_i(I_{n,i}) | \lambda(I_{n,i})\right] \big| > \frac{\epsilon R_n}{2 |F|} \right].
\end{align*}

\vspace{0.4 cm}

Now, by Chebyshev's inequality the left-hand side above is less than $$\frac{4 |F|^3}{\epsilon^2 R_n^2} \Big( \Var \big(N_i(I_{n,i})\big) + \Var\big( \E\left[N_i(I_{n,i}) | \lambda(I_{n,i})\right]\big) \Big),$$ and by the law of total variance $$\Var \big(N_i(I_{n,i})\big) = \E \left[\Var\big( N_i(I_{n,i}) | \lambda(I_{n,i})\big)\right] + \Var \big(\E\left[N_i(I_{n,i}) | \lambda(I_{n,i})\right]\big).$$

\vspace{0.4 cm}

But for the same reason as for (\ref{eq:espcondeqlambda}) we have $\Var \left(N_i(I_{n,i}) | \lambda(I_{n,i})\right) = \lambda(I_{n,i})$, so that $$\E \left[\Var\big( N_i(I_{n,i}) | \lambda(I_{n,i})\big)\right] \leq R_n.$$

\vspace{0.4 cm}

Therefore, it only remains to bound $$\frac{8 |F|^3}{\epsilon^2 R_n^2} \Var (\lambda_{i,n}).$$

\vspace{0.4 cm}

But then, by the same computations as in the proof of Lemma \ref{lemma:mainaux1} one gets
\begin{align*}
	\Var \left( \lambda_{i,n} \right) &= \Var \left(\int_t^{t + R_n} \xi_{n,i} (s) ds\right)\\
	&= \int_t^{t + R_n} \int_t^{t + R_n} \Cov \left(\xi_{n,i} (x), \xi_{n,i} (y)\right) dxdy,
\end{align*}

and then, using Lemma \ref{lemma:covfinite} $$\frac{8 |F|^3}{\epsilon^2 R_n^2} \Var (\lambda_{i,n}) \leq \frac{16 C_1 |F|^3}{C_2 \epsilon^2 R_n} + \frac{8 |F|^3}{\epsilon^2} \big( \P (\tau_n < t + R_n) + \epsilon_n\big),$$

and this last bound goes to $0$ as $n$ diverges.
\end{proof}

\section*{Acknowledgements}
This article was produced as part of the activities of FAPESP  Research, Innovation and Dissemination Center for Neuromathematics (grant 2013/ 07699-0). The author was supported by a FAPESP scholarship (grant 2020/12708-1).

\end{document}